\documentclass[11pt]{article}
\usepackage{graphicx}
\usepackage[a4paper]{geometry}
\usepackage[utf8]{inputenc}
\usepackage[T1]{fontenc}
\usepackage{aeguill}
\usepackage{amssymb}
\usepackage{amssymb,array}
\usepackage{syntax}
\usepackage{mdwtab}
\usepackage{mathenv}
\usepackage{float}
\usepackage{lscape}
\usepackage{here}
\usepackage[a4paper]{geometry}
\usepackage{amssymb}
\usepackage{bm}
\usepackage{fancyhdr}
\usepackage{graphicx}
\usepackage{natbib}
\usepackage{fancyvrb}

\newtheorem{theorem}{Theorem}
\newtheorem{lemma}[theorem]{Lemma}
\newtheorem{corollary}[theorem]{Corollary}

\newtheorem{remark}{Remark}

\newcommand{\cdf}[2]{\ensuremath{F_{#1}(#2)}}

\newcommand{\cdfEst}[3]{\ensuremath{	\widehat{F}_{#1} \left(#2 ; \Vect{#3} \right)}}
\newcommand{\cdfEstSimp}[1]{\ensuremath{	\widehat{F} \left(#1 \right)}}
\newcommand{\qtEstSimp}[1]{\ensuremath{	\widehat{\xi} \left(#1 \right)}}

\newcommand{\qt}[2]{\ensuremath{\xi_{#1}(#2)}}

\newcommand{\qtEst}[3]{\ensuremath{	\widehat{\xi}_{#1} \left(#2 ; \Vect{#3} \right)}}

\newcommand{\mypdf}[2]{\ensuremath{f_{#1}(#2)}}

\newcommand{\tauB}{\overline{\tau}}

\newcommand{\OO}[1]{\ensuremath{\mathcal{O}\left( #1 \right)}}
\newcommand{\Oas}[1]{\ensuremath{\mathcal{O}_{a.s.}\left( #1 \right)}}

\renewcommand{\Citet}[1]{\citeauthor{#1} \citeyearpar{#1}}

\newcommand{\Vect}[1]{\bm{#1}}
\newcommand{\RR}{\mathbb{R}}
\newcommand{\ZZ}{\mathbb{Z}}
\newcommand{\Esp}{\mathbf{E}}

\newcommand{\Prob}[1]{\mathbb{P}\left( #1\right)}
\newcommand{\NN}{\mathbb{N}}

\title{\Large \bf \sc Bahadur representation of sample quantiles for functional of Gaussian dependent sequences under a minimal assumption}

\author{ {\sc By Jean-Fran\c{c}ois Coeurjolly$^1$} \\
{\it SAGAG, Department of Statistics, Grenoble, FRANCE}}

\begin{document}
\maketitle

\begin{center}
{\small \begin{minipage}{12cm}
We obtain a Bahadur representation for sample quantiles of nonlinear functional of Gaussian sequences with correlation function decreasing as $k^{-\alpha}$ for some $\alpha>0$. This representation is derived under a mimimal assumption.

\end{minipage}}
\end{center}

\vspace*{.5cm}


\section{Introduction}

We consider the problem of obtaining a Bahadur representation of sample quantiles in a certain dependence context. Before stating in what a Bahadur representation consists, let us specify some general notation. Given some random variable $Y$, $\cdf{}{\cdot}=\cdf{Y}{\cdot}$ is referred as the cumulative distribution function of $Y$, $\qt{}{p}=\qt{Y}{p}$ for some $0<p<1$ as the quantile of order $p$. If \cdf{}{\cdot} is absolutely continuous with respect to Lebesgue measure, the probability density function is denoted by $\mypdf{}{\cdot}=\mypdf{Y}{\cdot}$. Based on the observation of a vector $\Vect{Y}=\left(Y(1),\ldots,Y(n) \right)$ of $n$ random variables distributed as $Y$, the sample cumulative distribution function and the sample quantile of order $p$ are respectively denoted by \cdfEst{Y}{\cdot}{Y} and \qtEst{Y}{p}{Y} or simply by \cdfEst{}{\cdot}{Y} and \qtEst{}{p}{Y}.

Let $\Vect{Y}=(Y(1),\ldots,Y(n))$ a vector of $n$ i.i.d. random variables such that $F^{\prime \prime}(\qt{}{p})$  exists and is bounded in a neighborhood of $\qt{}{p}$ and such that $F^\prime(\qt{}{p})>0$, Bahadur proved that as $n\to+\infty$,
$$
\qtEstSimp{p} - \qt{}{p} = \frac{p-\cdfEstSimp{p}}{\mypdf{}{\qt{}{p}}} + r_n,
$$
with $r_n=\Oas{n^{-3/4} \log(n)^{3/4}}$ where a sequence of random variables $U_n$ is said to be $\Oas{v_n}$ if $U_n/v_n$ is almost surely bounded. Kiefer obtained the exact rate  $n^{-3/4} \log \log(n)^{3/4}$. Under an Assumption on $F(\cdot)$ which is quite similar to the one done by Bahadur, extensions of above results to dependent random variables have been pursued in \Citet{Sen72} for $\phi-$mixing variables, in \Citet{Yoshihara95} for strongly mixing variables, and recently in \Citet{Wu05} for short-range and long-range dependent linear processes, following works of \Citet{Hesse90} and \Citet{Ho96}. Finally, such a representation has been obtained by \Citet{Coeurjolly07} for nonlinear functional of Gaussian sequences with correlation function decreasing as $k^{-\alpha}$ for some $\alpha>0$. 

\Citet{Ghosh71} proposed in the i.i.d. case a much simpler proof of Bahadur's result which suffices for many statistical applications. He established under a weaker assumption on $\cdf{}{\cdot}$ ($F^\prime(\cdot)$ exists and is bounded in a neighborhood of $\qt{}{p}$ and $\mypdf{}{\qt{}{p}}>0$) that the remainder term satisfies $r_n=o_{_{\mathbb{P}}}(n^{-1/2})$, which means that $n^{1/2}r_n$ tends to $0$ in probability. This result is sufficient for example to establish a central limit theorem for the sample quantile. Our goal is to extend Ghosh's result to nonlinear functional of Gaussian sequences with correlation function decreasing as $k^{-\alpha}$. The Bahadur representation is presented in Section~\ref{sec-result} and is applied to a central limit theorem for the sample quantile. Proofs are deferred in Section~\ref{sec-proofs}.

\section{Main result} \label{sec-result}

Let $\left\{ Y(i) \right\}_{i=1}^{+\infty}$ be a stationary (centered) gaussian process with variance 1, and correlation function $\rho(\cdot)$ such that, as $i\to +\infty$
\begin{equation} \label{hyporho}
| \rho(i) |\sim \;  i^{-\alpha}
\end{equation}
for some $\alpha>0$.

Let us recall some background on Hermite polynomials: the Hermite polynomials form an orthogonal system for the Gaussian measure and are in particular such that $\Esp\left( H_j(Y) H_k(Y) \right) = j!  \; \delta_{j,k}$, where $Y$ is referred to a standard Gaussian variable. For some measurable function $g(\cdot)$ defined on $\RR$ such that $\Esp(g(Y)^2)<+\infty$, the following expansion holds
$$
g(t) =\sum_{j \geq \tau} \frac{c_j}{j!} \; H_j(t) \quad \mbox{ with } \quad c_j = \Esp\left( g(Y) H_j(Y) \right),
$$
where the integer $\tau$ defined by $\tau= \inf \big\{ j \geq 0, c_j \neq 0 \big\}$,
is called the Hermite rank of the function~$g$. Note that this integer plays an important role. For example, it is related to the correlation of $g(Y_1)$ and $g(Y_2)$, for $Y_1$ and $Y_2$ two standard gaussian variables with correlation $\rho$, since $\Esp(g(Y_1)g(Y_2)=\sum_{k\geq \tau} \frac{(c_k)^2}{k!} \rho^k=\OO{\rho^\tau}$.

Our result is based on the assumption that $F_{g(Y)}^{ \prime}(\cdot)$ exists and is bounded in a neighborhood of $\qt{}{p}$. This is achieved if the function $g(\cdot)$ satisfies the following assumption (see {\it e.g.} \Citet{Dacunha82}, p.33).

$\bm{Assumption \; A(\qt{}{p}):}$ there exist $U_i$, $i=1,\ldots,L$, disjoint open sets such that $U_i$ contains a unique solution to the equation $g(t)=\qt{g(Y)}{p}$, such that $F^\prime_{g(Y)}(\qt{}{p})>0$ and such that $g$ is a $\mathcal{C}^1-$diffeomorphism on ${\displaystyle \cup_{i=1}^L U_i}$.

\noindent Note that this assumption allows us to obtain
$$
F_{g(Y)}^\prime( \qt{g(Y)}{p} )= \mypdf{g(Y)}{\qt{g(Y)}{p}} \; = \;\sum_{i=1}^L \frac{\phi(g_i^{-1}(t))}{g^\prime(g_i^{-1}(t))},
$$
where $g_i(\cdot)$ is the restriction of $g(\cdot)$ on $U_i$ and where $\phi(\cdot)$ is referred to the probability density function of a standard Gaussian variable.

Now, define, for some real $u$, the function $h_u(\cdot)$ by:
\begin{equation} \label{defh}
h_u(t) = \bm{1}_{\{  g(t) \leq u  \}}(t) - \cdf{g(Y)}{u}.
\end{equation}
We denote by $\tau(u)$ the Hermite rank of $h_{u}(\cdot)$. For the sake of simplicity, we set $\tau_p=\tau(\qt{g(Y)}{p})$. For some function $g(\cdot)$ satisfying Assumption $\bm{A(\qt{}{p})}$, we denote by
\begin{equation} \label{tauBarre}
\tauB_p = \inf_{\gamma \in \cup_{i=1}^L g(U_i)} \tau(\gamma),
\end{equation}
that is the minimal Hermite rank of  $h_u(\cdot)$ for $u$ in a neighborhood of $\qt{g(Y)}{p}$. Denote also by $c_j(u)$ the $j$-th Hermite coefficient of the function $h_u(\cdot)$.

\begin{theorem} \label{bahadur}
Under Assumption $\bm{A(\qt{}{p})}$, the following result holds as $n \to +\infty$
\begin{equation} \label{repBahadur}
\qtEst{}{p}{g(Y)}  - \qt{g(Y)}{p} = \frac{ p-\cdfEst{}{\qt{g(Y)}{p}}{g(Y)}   }{
\mypdf{g(Y)}{\qt{g(Y)}{p}}} \; + \; o_{_{\mathbb{P}}}\left(r_n(\alpha,\tauB_p) \right),
\end{equation}
where $\Vect{g(Y)}=\left(g(Y(1),\ldots,g(Y(n)) \right)$, for $i=1,\ldots, n$ and where the sequence $\big(r_n(\alpha,\tauB_p)\big)_{n\geq 1}$ is defined by
\begin{equation}\label{rn}
r_n( \alpha, \tauB_p) =\left\{
\begin{array}{ll}
n^{-1/2} & \mbox{ if } \alpha\tauB_p>1, \\
n^{-1/2}\log(n)^{1/2}         & \mbox{ if } \alpha\tauB_p=1, \\
n^{-\alpha\tauB_p/2}	 & \mbox{ if } \alpha\tauB_p<1.
\end{array}
\right.
\end{equation}
\end{theorem}

\begin{remark} \label{rem-rn}
The sequence $r_n(\alpha,\tauB_p)$ is related to the behaviour short-range or long-range dependent behaviour of the sequence $h_u(Y(1)),\ldots,h_u(Y(n))$ for $u$ in a neighborhood of $\qt{}{p}$. More precisely, it corresponds to the asymptotic behaviour of the sequence
$$
\left( \frac1n \sum_{|i|<n} \rho(i)^{\tauB_p} \right)^{1/2}.
$$
\end{remark}

\begin{corollary}
Under Assumption $\bm{A(\qt{}{p})}$, then the following convergence in distribution hold as $n \to +\infty$

$(i)$ if $\alpha \tauB_p>1$
\begin{equation} \label{convLoiQt1}
\sqrt{n} \left( \qtEst{}{p}{g(Y)}  - \qt{g(Y)}{p} \right) \stackrel{d}{\longrightarrow} \mathcal{N}(0,\sigma^2_p),
\end{equation}
where 
$$\sigma^2_p=\frac1{f(p)^2}\sum_{i \in \ZZ} \sum_{j \geq \tauB_p} \frac{c_j(p)^2}{j!} \rho(i)^j\mbox{ with } f(p)=f_{g(Y)}(\qt{g(Y)}{p}) \mbox{ and } c_j(p)=c_j(\qt{g(Y)}{p}).$$

$(ii)$ if $\alpha\tauB_p<1$
\begin{equation} \label{convLoiQt2}
n^{\alpha \tauB_p/2} \left( \qtEst{}{p}{g(Y)}  - \qt{g(Y)}{p} \right) \stackrel{d}{\longrightarrow} \frac{c_{\tauB_p}(p)}{\tauB_p!f(p)} Z_{\tauB_p},
\end{equation}
where
$$
Z_{\tauB_p} = K(\tauB_p,\alpha) \int^{\prime}_{\RR^{\tauB_p}} \frac{\exp(i(\lambda_1+\cdots+\lambda_{\tauB_p}))-1}{i(\lambda_1+\cdots+\lambda_{\tauB_p})} \prod_{j=1}^{\tauB_p} |\lambda_j|^{(\alpha-1)/2}\widetilde{B}(d\lambda_j)
$$
and
$$
K(\tauB_p,\alpha)= \left( 
\frac{(1-\alpha\tauB_p/2)(1-\alpha\tauB_p)}{\tauB_p! 
\left( 2\Gamma(\alpha) \sin(\pi(1-\alpha)/2) \right)^{\tauB_p} }
\right)^{1/2}.
$$
The measure $\widetilde{B}$ is a Gaussian complex measure and the symbol $\int^\prime$ means that the domain of integration excludes the hyperdiagonals $\{\lambda_i= \pm \lambda_j, i\neq j\}$.
 
\end{corollary}

The proof of this result is omitted since it is a direct application of Theorem~\ref{bahadur} and general limit theorems adapted to nonlinear functional of Gaussian sequences, {\it e.g.} \Citet{Breuer83} and \Citet{Dehling89}.

\section{Proofs} \label{sec-proofs}

\subsection{Auxiliary Lemma}

\begin{lemma} \label{calculInt}
For every $j\geq 1$ and for all positive sequence $(u_n)_{n\geq 1}$ such that $u_n\to 0$, as $n \to +\infty$, we have, under Assumption $\mathbf{A(\qt{}{p})}$
\begin{equation} \label{I}
I= \int_{\RR} H_j(t) \phi(t) \mathbf{1}_{ \{  |g(t)-\qt{g(Y)}{p})|\leq u_n  \}} dt \sim  u_n \; \kappa_j ,
\end{equation}
where $\kappa_j$ is defined, for every $j\geq 1$,by
\begin{equation} \label{kappaj}
\kappa_j = \left\{
\begin{array}{ll}
-2  \sum_{i=1}^L \frac{\phi^{\prime}\left( g_i^{-1}(\qt{}{p}\right)}{g^\prime \left(  g_i^{-1}(\qt{}{p})   \right)} & \mbox{ if } j=1 ,\\
2 (-1)^j  \sum_{i=1}^L \frac{\phi^{(j-2)}\left( g_i^{-1}(\qt{}{p}\right)}{g^\prime \left(  g_i^{-1}(\qt{}{p})   \right)} & \mbox{ if } j >1.\\
\end{array} \right.
\end{equation}
\end{lemma}

\begin{proof}
Under Assumption $\mathbf{A(\qt{}{p})}$, there exists $n_0 \in \NN$ such that for all $n\geq n_0$,
\begin{equation} \label{eq1I}
I = \sum_{i=1}^L I_i
\quad \mbox{ with } \quad
I_i = \int_{U_i} H_j(t) \phi(t) \mathbf{1}_{ \{  \qt{}{p}-u_n \; \leq\;  g(t)\; \leq \;\qt{}{p}+u_n  \}} dt.
\end{equation}
Assume without loss of generality that the restriction of $g(\cdot)$ on $U_i$ (denoted by $g_i(\cdot)$) is an increasing function, we have
\begin{eqnarray}
I_i &=&\int_{U_i} H_j(t) \phi(t) \mathbf{1}_{ \{  \qt{}{p}-u_n \; \leq\;  g(t)\; \leq \;\qt{}{p}+u_n  \}} dt \nonumber \\&=&
\int_{g_i^{-1}(\qt{}{p}-u_n)}^{g_i^{-1}(\qt{}{p}+u_n)} H_j(t) \phi(t) dt \nonumber \\
&=& \left\{ \begin{array}{ll}
\phi(m_{i,n})- \phi(M_{i,n})  = (m_{i,n}-M_{i,n}) &  \mbox{ if } j=1 \\
(-1)^j \left(  \phi^{(j-1)}(  M_{i,n}     ) - \phi^{(j-1)}( m_{i,n}    ) \right)&  \mbox{ if } j>1,
\end{array} \right. \nonumber
\end{eqnarray}
where $M_{i,n}=g_i^{-1}(\qt{}{p}+u_n)$ and $m_{i,n}=g_i^{-1}(\qt{}{p}-u_n)$. Then, there exists $\omega_{n,i,j} \in [m_{i,n},M_{i,n}]$ for every $j\geq 1$ such that
$$
I_i = \left\{ \begin{array}{ll}
(m_{i,n}-M_{i,n}) \; \phi^{(1)}(\omega_{n,i,1}) & \mbox{ if } j=1 \\
 (-1)^j \left( M_{i,n}-m_{i,n} \right) \phi^{(j-2)} (\omega_{n,i,j}) &  \mbox{ if } j>1.
\end{array} \right. ,
$$
Under Assumption $\mathbf{A(\qt{}{p})}$, we have, as $n \to +\infty$
$$
 \omega_{n,i,j} \sim g_i^{-1}(\qt{}{p}) \quad \mbox{ and } \quad
M_{i,n} - m_{i,n} \sim  2 u_n \; \frac{1}{g^\prime(g_i^{-1}(\qt{}{p})) },
$$
which ends the proof.
\end{proof}

\subsection{Proof of Theorem~\ref{bahadur}}

For the sake of simplicity, we set $\qtEstSimp{p}=\qtEst{}{p}{g(Y)}$, $\qt{}{p}=\qt{g(Y)}{p}$, $\cdfEstSimp{\cdot}=\cdfEst{}{\cdot}{g(Y)}$, $\cdf{}{\cdot}=\cdf{g(Y)}{\cdot}$ et $\mypdf{}{\cdot}=\mypdf{g(Y)}{\cdot}$ and $r_n=r_n(\alpha,\tauB_p)$. Define,
$$
V_n= r_n^{-1} \left( \qtEstSimp{p} - \qt{}{p} \right) \quad \mbox{and} \quad W_n= r_n^{-1} \left( \frac{p-\cdf{}{p}}{\mypdf{}{p}} \right).
$$
The result is established if $V_n-W_n \stackrel{\mathbb{P}}{\rightarrow} 0$ as $n \to+\infty$. It suffices to prove that $V_n$ and $W_n$ satisfy the conditions of Lemma~1 of \Citet{Ghosh71}:
\begin{list}{$\bullet$}{}
\item \textbf{condition (a)}~: for all $\delta>0$, there exists $\varepsilon=\varepsilon(\delta)$ such that $\Prob{|W_n|>\varepsilon}<~\delta$.
\item \textbf{condition (b)}~: for all $y \in \RR$ and for all $\varepsilon>0$
$$\lim_{n\to+\infty} \Prob{V_n \leq y, W_n \geq k+\varepsilon} \quad \mbox{ and } \quad
\lim_{n\to+\infty} \Prob{V_n \geq y+\varepsilon, W_n \geq k} $$
\end{list}

\noindent \textbf{condition (a)}~: from Bienaym\'e-Tchebyshev's inequality it is sufficient to prove that $\Esp W_n^2 = \mathcal{O}(1)$. Rewrite $W_n = \frac{r_n^{-1}}{n} \sum_{i=\ell+1}^n h_{\qt{}{p}}\left( Y(i)   \right)$. Let $c_j$ (for some $j\geq 0$) denote the $j$-th Hermite coefficient of $h_{\qt{}{p}}(\cdot)$. Since $h_{\qt{}{p}}(\cdot)$ has at least Hermite rank $\tauB_p$, then
\begin{eqnarray}
\Esp W_n^2 &=& \frac{r_n^{-2} }{n^2} \sum_{i_1,i_2=1}^n \Esp\left( h_{\qt{}{p}}\left( Y(i_1) \right)  h_{\qt{}{p}}\left( Y(i_2) \right)      \right) \nonumber \\
&=& \frac{r_n^{-2} }{n^2} \sum_{i_1,i_2=1}^n  \sum_{j_1,j_2\geq \tau_p} c_{j_1} c_{j_2}  \Esp\left(
H_{j_1} \left( Y(i_1) \right) H_{j_2} \left( Y(i_2) \right) \right) \nonumber \\
&=& \frac{r_n^{-2} }{n^2} \sum_{i_1,i_2=1}^n \sum_{j\geq \tauB_p} \frac{(c_{j})^2}{ (j)!} \; \rho(i_2-i_1)^{j} \nonumber \\
&=& \mathcal{O} \left(  r_n^{-2}  \times \frac1n \sum_{|i|<n} \rho(i)^{\tauB_p}   \right) = \OO{1},\nonumber
\end{eqnarray}
from Remark~\ref{rem-rn}.

\noindent \textbf{condition (b)}~: let $y \in \RR$, we have
\begin{eqnarray}
\left\{ V_n \leq y \right\} &=& \left\{ \qtEstSimp{p} \leq y \times r_n  + \qt{}{p}    \right\} \nonumber \\
&=& \left\{ p \leq \cdfEstSimp{y\times r_n + \qt{}{p}} \right\} = \left\{ Z_n \leq y_n \right\} ,
\label{condbEq1}
\end{eqnarray}
with
$$
Z_n = \frac{r_n^{-1} }{\mypdf{}{\qt{}{p}}} \left( F\Big( y\times r_n + \xi(p)\Big) - \cdfEstSimp{\frac{y}{\sqrt{r_n }} + \xi(p)} \right)$$
and $$ y_n= \frac{r_n^{-1}}{f(\xi(p))}  \Big(  F\Big( y\times r_n  + \xi(p)\Big)-p  \Big)
$$
Under Assumption $\mathbf{A(\qt{}{p})}$, we have $y_n\to y$, as $n \to +\infty$. Now, prove that $Z_n-W_n \stackrel{\mathbb{P}}{\rightarrow} 0$. Without loss of generality, assume $y>0$. Then, we have
\begin{eqnarray}
W_n-Z_n &=& \frac{r_n^{-1}}{\mypdf{}{p}} \left( \widehat{F}\Big( y\times r_n  + \qt{}{p}\Big) - F\Big(y\times r_n  + \qt{}{p}\Big) - \cdfEstSimp{ \qt{}{p}} + \cdf{}{\qt{}{p} }  \right)  \nonumber \\
&=& \frac{r_n^{-1}}{n} \; \frac{1}{\mypdf{}{\qt{}{p}}} \; \sum_{i=1}^n h_{\qt{}{p},n} \left( Y(i)\right) \nonumber
\end{eqnarray}
where $h_{\qt{}{p},n}(\cdot)$ is the function defined for $t\in \RR$ by~:
$$
h_{\qt{}{p},n}(t) = \mathbf{1}_{\Big\{ \qt{}{p}\leq g(t) \leq \qt{}{p}+y\times r_n  \Big\} } (t) - \mathbb{P} \Big( \qt{}{p}\leq g(Y)\leq \qt{}{p}+  y\times r_n  \Big).
$$
For $n$ sufficiently large, the function $h_{\qt{}{p},n}(\cdot)$ has Hermite rank $\tauB_p$. Denote by $c_{j,n}$ the $j$-th Hermite coefficient of $h_{\qt{}{p},n}(\cdot)$. From Lemma~\ref{calculInt}, there exists a sequence $(\kappa_j)_{j\geq \tauB_p}$ such that, as $n \to +\infty$
$$
c_{j,n} \; \sim \; \kappa_j \times r_n.$$
Since, for all $n\geq 1$ $\Esp(h_n(Y)^2) = \sum_{j \geq \tauB_p} (c_{j,n})^2 / j!< +\infty$, it is clear that the sequence $(\kappa_j)_{j\geq \tauB_p}$ is such that $\sum_{j \geq \tauB_p} (\kappa_j)^2/ j! <+\infty$. By denoting $\lambda$ a positive constant, we get, as $n \to +\infty$
\begin{eqnarray}
\Esp (W_n-Z_n)^2 &=& \frac{r_n^{-2} }{n^2} \; \frac{1}{ \mypdf{}{\qt{}{p}}^2} \sum_{i_1,i_2=1}^n \Esp \left( h_{\qt{}{p},n}\left( Y(i_1)\right) h_{\qt{}{p},n}\left( Y(i_2)\right) \right) \nonumber \\
&=& \frac{r_n^{-2} }{n^2} \; \frac{1}{ \mypdf{}{\qt{}{p}}^2} \sum_{i_1,i_2=1}^n \sum_{j_1,j_2\geq \tauB_p} c_{j_1,n} c_{j_2,n}
\Esp\left( H_{j_1}\left(Y(i_1)\right)H_{j_2}\left(Y(i_2)\right) \right) \nonumber \\
&=& \frac{r_n^{-2} }{n^2} \; \frac{1}{\mypdf{}{\qt{}{p}}^2} \sum_{i_1,i_2=1}^n \sum_{j \geq \tauB_p} \frac{c_{j,n}^2}{ j!} \rho(i_2-i_1)^{j} \nonumber \\
&\leq & \lambda \; \frac{r_n^{-2} }{n} \sum_{j \geq \tauB_p} \frac{(\kappa_{j})^2}{ j! } \; r_n^2 \sum_{|i|<n} \rho(i)^{j} = \mathcal{O}\left(\frac1n \sum_{|i|<n} \rho(i)^{\tauB_p} \right) = \mathcal{O}(r_n^2), \nonumber
\end{eqnarray}
from Remark~\ref{rem-rn}. Therefore, $W_n-Z_n$ converges to 0 in probability, as $n\to+\infty$. Thus, for all $\varepsilon>0$, we have, as $n \to +\infty$,
$$
\Prob{ V_n\leq y , W_n\geq y+\varepsilon} = \Prob{Z_n \leq y_n , W_n \geq y+\varepsilon   } \rightarrow 0.
$$
Following the sketch of this proof, we also have $\Prob{ V_n\geq y+\varepsilon , W_n\leq y} \to 0$, ensuring condition~(b).
Therefore, $W_n-Z_n$ converges to 0 in probability, as $n\to+\infty$. Thus, for all $\varepsilon>0$, we have, as $n \to +\infty$,
$$
\Prob{ V_n\leq y , W_n\geq y+\varepsilon} = \Prob{Z_n \leq y_n , W_n \geq y+\varepsilon   } \rightarrow 0.
$$
Following the sketch of this proof, we also have $\Prob{ V_n\geq y+\varepsilon , W_n\leq y} \to 0$, ensuring condition~(b).

\bigskip

\noindent{\textbf{\large Acknowledgement}} The author would like to thank the referee for his coments improving the statement of Corollary~2.

\bigskip \bigskip

\noindent {\sc J.-F. Coeurjolly \\
SAGAG, Department of Statistics, Grenoble \\
1251 Av. Centrale  BP 47\\
38040 GRENOBLE Cedex 09\\
France \\
E-mail:} Jean-Francois.Coeurjolly@upmf-grenoble.fr

\end{document}